\documentclass[11pt]
{amsart}
\usepackage[margin=1.4in]{geometry}

\usepackage{amsmath, amssymb, amsfonts, amsthm, enumerate}

\newtheorem{theorem}{\bf Theorem}[section]

\newtheorem{conjecture}[theorem]{\bf Conjecture}

\newtheorem{proposition}[theorem]{\bf Proposition}

\let\theta=\vartheta
\let\rho=\varrho
\let\sigma=\varsigma

\begin{document}

\title[Discrepancy of high-dimensional permutations]{Discrepancy of high-dimensional permutations}

\author[N.~Linial]{Nati Linial}
\address{School of Computer Science and engineering, The Hebrew University, Jerusalem, Israel.}
\thanks{Supported by ERC grant 339096 High-Dimensional Combinatorics}
\email{nati@cs.huji.ac.il}

\author[Z.~Luria]{Zur Luria}
\address{Institute of Theoretical Studies, ETH, 8092 Zurich, Switzerland.}
\email{zluria@gmail.com}

\date{\today}

\begin{abstract}
Let $L$ be an order-$n$ Latin square. For $X, Y, Z \subseteq \{1, ... ,n\}$, let $L(X, Y, Z)$ be the number of triples $i\in X, j\in Y, k\in Z$ such that $L(i,j) = k$. We conjecture that asymptotically almost every Latin square satisfies $|L(X, Y, Z) - \frac 1n |X||Y||Z||= O(\sqrt{|X||Y||Z|})$ for every $X, Y$ and $Z$. Let $\varepsilon(L):= \max |X||Y||Z|$ when $L(X, Y, Z)=0$. The above conjecture implies that $\varepsilon(L) = O(n^2)$ holds asymptotically almost surely (this bound is obviously tight). We show that there exist Latin squares with $\varepsilon(L) = O(n^2)$, and that $\varepsilon(L) = O(n^2 \log^2 n)$ for almost every order-$n$ Latin square. On the other hand, we recall that $\varepsilon(L)\geq \Omega(n^{33/14})$ if $L$ is the multiplication table of an order-$n$ group. We also show the existence of Latin squares in which every empty {\em cube} has side length $O((n\log n)^{1/2})$, which is tight up to the $\sqrt{\log n}$ factor. Some of these results extend to higher dimensions. Many open problems remain.
\end{abstract}

\maketitle

\section{Introduction}

The notion of {\em discrepancy} is central to all branches of discrete mathematics. Indeed, several books~\cite{Mat, Chaz, BeSo} have been dedicated to this subject. Roughly speaking, one asks how well finite sets can approximate a uniform measure. A bit more concretely, the problem is defined in terms of a collection $\mathcal F$ of subsets in a probability space $(\Omega,\mu)$. We seek the minimum of $\sup_{X\in \mathcal F} |\frac{|S\cap X|}{|S|}-\mu(X)|$ over all sets $S$ of given cardinality. Such questions and their many variants make sense and are interesting in numerous contexts. An important example from graph theory is the {\em expander mixing lemma}. Let $G=(V,E)$ be a $d$-regular $n$-vertex graph. This lemma asserts that if $G$ is an expander graph, then for every two subsets $A, B\subseteq V$ there holds $|e(A,B)-\frac dn|A||B||= O(\sqrt{|A||B|})$ where $e(A,B)$ is the number of ordered pairs $(a, b)$ with $a\in A, b\in B$ and $ab\in E$. The unspecified constant in the big-oh term depends on $G$'s spectrum, but we do not elaborate on this point and refer the reader to the survey~\cite{hlw}.

A considerable body of recent research is aimed at developing a theory of {\em high-dimensional combinatorics}. Many basic combinatorial constructs have interesting high-dimensional counterparts, and it is natural to study discrepancy phenomena in these frameworks. Specifically we consider discrepancy in {\em high-dimensional permutations}. Let us briefly recall this concept~\cite{LiLu14}. We equate a (classical, i.e., one-dimensional) permutation with its permutation matrix, namely, an $n\times n$ array of zeros and ones where every row and every column contains exactly one $1$. In analogy, a $d$-dimensional permutation $A$ is an $[n]^{d+1}=n\times n\times\ldots \times n$ array of zeros and ones such that for every index $d+1\ge i\ge 1$ and every choice of integers $\alpha_j\in [n]$ over $1\le j\neq i \le d+1$ there is exactly one choice of $x\in [n]$ for which $A(\alpha_1,\ldots, \alpha_{i-1}, x, \alpha_{i+1},\ldots, \alpha_{d+1})=1$. Note, in particular, that a two-dimensional permutation is synonymous with a Latin square.

The class $\mathcal F$ that defines our discrepancy problem is comprised of all boxes $\mathcal{T}=T_1\times\ldots\times T_{d+1}\subseteq [n]^{d+1}$. The {\em volume} of this box is defined to be $\text{vol}(\mathcal{T}):=\prod |T_i|$. Our discrepancy problem is to find $d$-dimensional permutations $A$, such that for every box $\mathcal T$ it holds that $A(\mathcal{T}):=|\{\alpha \in \mathcal{T}:A(\alpha)=1\}|$ is close to $\frac{\text{vol}(\mathcal{T})}{n}$. (Clearly this is what one would expect, since the density of $1$ entries in a $d$-dimensional permutation is $\frac 1n$).
We propose the following conjecture.

\begin{conjecture}
\label{conj:high_dim_EML}
For every $d\geq 2$ there exist arbitrarily large $d$-dimensional permutations $A$ such that for every box  $\mathcal T$ we have
\[
\left|A(\mathcal{T}) - \frac{\text{vol}(\mathcal{T})}{n}\right| = O(\sqrt{\text{vol}(\mathcal T)}).
\]
\end{conjecture}

There are at least two reasons why we expect this to be true. Consider the following ``poor man's analog" of a random Latin square. It is a random $n\times n\times n$ array of zeros and ones whose entries are chosen independently with the same distribution, where $1$ is chosen with probability $\frac{1}{n}$. It is easily verified that this relation holds in that model. In addition, a $d$-dimensional permutation may be viewed as a $(d+1)$-partite $(d+1)$-uniform hypergraph, and we find the similarity with the expander mixing lemma rather compelling.

We say that $\mathcal{T}$ is an \textit{empty box} in $A$ if $A(\mathcal{T}) = 0$, and denote by $\varepsilon(A)$ the maximal volume of an empty box in $A$. One consequence of the above conjecture is that there are $d$-dimensional permutations $A$ such that $\varepsilon(A) = O(n^2)$. On the other hand, it is easy to see that $\varepsilon(A) = \Omega(n^2)$ for {\em every} $d$-dimensional permutation, since every (classical) permutation matrix contains a $\lfloor\frac n2\rfloor\times\lfloor\frac n2\rfloor$ block of zeros. Indeed, let  $A$ be an arbitrary $d$-dimensional permutation. Pick some $T_2\subseteq [n]$ of cardinality $\lfloor\frac n2\rfloor$ and some $t_3,\dots, t_{d+1}\in [n]$, and let $T_3=\{t_3\},\ldots, T_{d+1}=\{t_{d+1}\}$. We can find a subset $T_1\subseteq [n]$ of cardinality $\lfloor\frac n2\rfloor$ for which $\mathcal{T}=T_1\times\ldots\times T_{d+1}\subseteq [n]^{d+1}$ is an empty box in $A$. Indeed, for every $t\in T_2$, there is exactly one $x\in [n]$ for which $A(x, t, t_3,\ldots, t_{d+1})=1$ and clearly $x$ cannot belong to $T_1$. But altogether only $\lfloor\frac n2\rfloor$ elements are ruled out from being in $T_1$, one per each element of $T_2$ so that at least $\lfloor\frac n2\rfloor$ are still acceptable and the claim follows.

We prove the following theorems in this spirit for $2$-dimensional permutations, i.e., for Latin squares.
\begin{theorem}
\label{thm:typical}
Asymptotically almost every order-$n$ Latin square $A$ satisfies $\varepsilon(A) = O(n^2 \log^2(n))$.
\end{theorem}

\begin{theorem}
\label{thm:upper_bound}
There exist infinitely many order-$n$ Latin squares satisfying $\varepsilon(A) = O(n^2)$ (and hence $\varepsilon(A) = \Theta(n^2)$).
\end{theorem}

We tend to believe the following statement which subsumes both theorems:

\begin{conjecture}
Asymptotically almost every order-$n$ Latin square $A$ satisfies $\varepsilon(A) = O(n^2)$.
\end{conjecture}

Moreover, it is conceivable that the discrepancy condition of Conjecture~\ref{conj:high_dim_EML} holds for asymptotically almost every $d$-dimensional permutation.

It is easy to see that the multiplication table of a finite group is a Latin square, and  problems that we consider here have been previously addressed in the group theory literature. Babai and Sos~\cite{BaSo85}, defined a subset $S\subset\Gamma$ of a finite group to be {\em product-free} if there are no three elements $x,y ,z\in S$ with $xy=z$. Note that in our language this means that $S\times S\times S$ is an empty box in the Latin square $L$ corresponding to $\Gamma$. Using the classification of finite simple groups, Babai and Sos showed that every finite group contains large product-free sets. Subsequently, Kedlaya~\cite{Ke97} improved their bound. His result implies:

\begin{theorem}[Kedlaya]
If $L$ is a Latin square that is the mutiplication table of an order-$n$ group, then $\varepsilon(L) \ge c n^{\frac{33}{14}}$ for some fixed $c>0$.
\end{theorem}

On the other hand, Gowers~\cite{Go08} has exhibited order-$n$ groups for which $\varepsilon(L) \le C n^{\frac{8}{3}}$ for some fixed $C>0$.

These results show that a typical Latin square has substantially lower discrepancy than any group of the same order.

A {\em cube} is a box $A\times B\times C$ with $|A|=|B|=|C|$. It is easy to see that every order-$n$ Latin square has an empty cube of side $\lfloor(n+1/4)^{1/2} -1/2\rfloor$, and we can show the following.

\begin{theorem}\label{empty:cube}
There exist infinitely many order-$n$ Latin squares $L$ in which every empty cube has side $O((n \log n)^{1/2})$.
\end{theorem}

As mentioned, Kedlaya finds an empty cube of side $\Omega(n^{11/14})$ in the Latin square of every order-$n$ group.

Again, analogs in general dimension suggest themselves as we discuss in Section~\ref{section:open}.

The proof of Theorem~\ref{thm:typical} is based on our earlier work~\cite{LiLu14} in which we derived an upper bound on the number of $d$-dimensional permutations. The proof of Theorems~\ref{thm:upper_bound} and~\ref{empty:cube} is based on ideas developed by P. Keevash in his recent breakthrough work on the theory of combinatorial designs. He considers in~\cite{Ke15} a random greedy process in which a set system evolves as sets are added to it in sequence. As he shows, with high probability the partial design that is obtained this way can be completed to a bona-fide design.

\section{Proof of theorem \ref{thm:typical}}

This result follows from an upper bound on high dimensional permanents proved in \cite{LiLu14}. Recall that the support of an $r$-dimensional array $X$ is 
\[
\text{Supp}(X) = \{(i_1, ... ,i_r):X(i_1, ... ,i_r)\ne 0\}.
\]
Define the permanent of a $(d+1)$-dimensional 0-1 array $A$ to be the number of $d$-permutations whose support is contained in $\text{Supp}(A)$, and let $r_{i_1,...,i_d}$ be the number of ones in the line $A(i_1,...,i_d,\cdot)$, i.e., the number of integers $x\in[n]$ for which $A(i_1,...,i_d,x)=1$. Then
\begin{equation}
\label{permanent_bound}
\text{Per}(A) \leq \prod_{i_1,...,i_d=1}^n{\left(1+O\left(\frac{\log^d(r_{i_1,...,i_d})}{r_{i_1,...,i_d}}\right)\right)\frac{r_{i_1,...,i_d}}{e^d}} .
\end{equation}

We denote the number of order-$n$ Latin squares by $\mathcal{L}(n)$.
Fix sets $X, Y, Z \subseteq [n]$ and let $B$ denote the $n \times n \times n$ 0-1 array which is $0$ in the box $X \times Y \times Z$ and $1$ elsewhere. The probability that $X  \times Y \times Z$ is an empty box of an order-$n$ Latin square chosen uniformly at random is $\frac{\text{Per}(B)}{\mathcal{L}(n)}$.
A counting argument due to van Lint and Wilson \cite{VLW} shows that $\mathcal{L}(n) = \left((1+o(1))\frac{n}{e^2}\right)^{n^2}$ and in particular $\mathcal{L}(n)\geq \left( \frac{n}{e^2}\right)^{n^2}$, and so we obtain the following upper bound by applying (\ref{permanent_bound}) to $B$.
\[
\Pr(L(X,Y,Z)=0) \leq (1+O(\log^2 n/n))^{n^2} \cdot
\frac{\left(\frac{n}{e^2}\right)^{n^2-|X||Y|}\cdot \left(\frac{n-|Z|}{e^2}\right)^{|X||Y|}}{ \left( \frac{n}{e^2}\right)^{n^2}}
\]
\[
\leq e^{O(n\log^2(n))} e^{-|X||Y||Z|/n}.
\]

Next, we apply the union bound over all boxes whose volume is at least  $M n^2 \log^2(n)$ for a large constant $M$ whose value will be chosen later. There are at most $(2^n)^3$ ways to choose $A,B$ and $C$, and so if $L$ is an order-$n$ Latin square chosen uniformly at random, we have
\[
\Pr(\varepsilon(L)\geq M n^2 \log^2(n)) \leq 2^{3n} \cdot e^{O(n\log^2(n))} e^{-M n \log^2(n)}.
\]
Therefore, for any constant $M$ that is larger than the constant in the big-oh term, we obtain a vanishingly small probability.

\section{Proof of theorem \ref{thm:upper_bound}}

Here we use an insight from Keevash's breakthrough papers \cite{Ke14,Ke15} on the existence and asymptotic enumeration of designs. We consider his construction for the specific case of Steiner triple systems. The first part of the algorithm involves a random greedy strategy which is stopped when all but a vanishingly small fraction of the vertex pairs are covered by triples. The crux of his proof is that, with high probability, the resulting set of uncovered triples can be completed to a Steiner triple system.

An analogous result is most likely also true for the random construction of Latin squares. However, to simplify matters, we use Keevash's results on Steiner triple systems and adapt them to our needs. Note that every order-$n$ Steiner triple system $X$ yields a (symmetric) order-$n$ Latin square $L$ as follows: $L(i,j)=k \Leftrightarrow \{i,j,k\} \in X$ and $L(i,i)=i$ for all $i\in [n]$. We define an {\em empty box} in $X$ to be a triple of sets $A,B,C \subseteq [n]$ such that $\{i,j,k\}\not\in X$ for every $i \in A, j \in B, k \in C$. We say that this box has {\em volume} $|A||B||C|$, and denote the largest volume of an empty box in $X$ by $\varphi(X)$. Since an empty box in $L$ is also an empty box in $X$, we have $\varepsilon(L) \le \varphi(X)$.

Steiner triple systems constructed using Keevash's method tend to have no large empty boxes:

\begin{proposition}
\label{prop:keevash_empty}
Almost every order-$n$ Steiner triple system $X$ generated by Keevash's method satisfies $\varphi(X)\leq M n^2$. Here $M>0$ is an absolute constant, e.g., $M=9000$ will do.
\end{proposition}

Keevash's algorithm asymptotically almost surely constructs a Steiner triple system for every large enough $n$ such that $n \equiv 1$ or $3$ (mod $6$). This proposition implies that for such $n$ there exist order-$n$ Latin squares $L$ with $\varepsilon(L) \leq M n^2$.

\begin{proof}[Proof of the proposition: ]

In view of the way in which Keevash's construction proceeds, it suffices to show that at the end of the random greedy process there remain no large empty boxes. Since that process is monotone and triples only get added, it suffices to show that after a small fraction of this stage is completed, no large box remains empty. Recall that at each step of the process a triple is chosen at random from among the {\em legal} triples, i.e., those that have at most one vertex in common with every previously selected triple.

Given $A,B,C \subseteq [n]$, an $ABC$ triple is a triple that meets $A, B$ and $C$. An $AB\bar{C}$ triple meets $A$ and $B$, but does not meet $C$, etc. Clearly, the set $F$ of all $ABC$ triples satisfies $\frac{1-o(1)}6 |A||B||C|\leq|F|\leq |A||B||C|$. Let $\lambda>0$ be a constant whose value will be chosen later. We refer to the initial  $\lambda n^2$ steps of the random greedy process as the {\em first stage}, and prove that if $|A||B||C| \ge Mn^2$, then it is very unlikely that no triple in $F$ is selected during the first stage. There are $8^n$ choices for $A, B, C$, so if for every choice of $A, B, C$ this statement fails with probability $o(8^{-n})$, our claim is established.

Indeed, this sounds plausible. Since
$|F|/{n\choose 3} \geq (1-o(1))(M/n)$
, the probability that during $\lambda n^2$ steps we never select a triple from $F$ ought to be exponentially small in $n$.
However, this heuristic argument ignores the fact that triples in $F$ may become illegal during the process even if we never select a triple from $F$. Thus, the choice of an $AB\bar{C}$ triple may invalidate as many as $3|C|$ triples in $F$. We need to show that whp not too many such choices are made\footnote{We say informally that an event holds {\em with high probability} (whp) meaning that it holds with probability $\ge 1-p^n$ without specifying $p$. The relevant range of $p$ is given in our formal discussion.}.


We will show that
\begin{equation}\label{f_safe}
\text{Whp, at the end of first stage, at least~} \frac{|F|}{2} \text{~triples in~} F \text{~remain legal}.
\end{equation}
Consequently, the probability that during the first stage we select no member of $F$ is at most
\[
\left(1-\frac{|F|}{2{n\choose 3}}\right)^{\lambda n^2}=e^{-(1+o(1)) 3 \lambda |F|/n}.
\]
As $|F|\geq M n^2/6$, this is $o(8^{-n})$, provided that $\lambda M > 6 \ln 2$.

To prove Statement~(\ref{f_safe}) we show first that whp during the first stage at most $|A||B|/108$ triples of type $AB\bar{C}$ get chosen. Thus the chosen $AB\bar{C}$ triples invalidate at most $3|C|\cdot |A||B|/108\le |F|/6$ triples of $F$. Together with the analogous contribution of types $A\bar{B}C$ and $\bar{A}BC$, at most half of the triples in $F$ get invalidated, and so at least half of them remain legal.

There are at most $|A||B|n$ triples of type $AB\bar{C}$, and each chosen triple invalidates at most $3n$ triples. If $X$ is the number of type $AB\bar{C}$ triples that we sample during the first stage, then $X = \sum_{i=1}^{\lambda n^2}{X_i}$, where $X_i$ is the indicator random variable of the event that the $i$-th chosen triangle is in $AB\bar{C}$. Therefore,

\[\mathbb{E}X  = \sum_{i=1}^{\lambda n^2} \mathbb{E}X_i \le \frac{\lambda |A||B|n^3}{{n\choose 3}-3\lambda n^3}=\frac{6\lambda |A||B|}{1 -18\lambda-o(1)}.\]

We recall the following generalization of Chernoff's inequality (see Theorem 3.4 in \cite{PS97}).
Namely, if $Y$ is the sum of $N$ Bernoulli random variables $Y_1, ... ,Y_N$ and for every subset $S \subset [N]$ we have $\Pr(\wedge_{i \in S}Y_i=1) \leq p^{|S|}$ for some $0<p<1$, then for every $\delta>0$, 
\[
\Pr(Y\geq (1+\delta) N p) \leq \exp\left(-\frac{\delta^2}{2+\delta} N p\right).
\]

Let $q := |A||B|n/\left(\binom{n}{3}-3 \lambda n^3\right)$.
The probability that $X_i = 1$ conditioned on the values of previous variables is always at most $q$, and so $\Pr(\wedge_{i \in S}X_i=1) \leq q^{|S|}$ for every $S \subset \lambda n^2$. Moreover, $|A||B|/108\ge K\cdot\mathbb{E}X$, where $K=\frac{1-18\lambda-o(1)}{648\lambda}$, and so we have
\[\Pr(X>|A||B|/108)\le \exp\left(-\frac{(K-1)^2}{K+1}\cdot\frac{6\lambda |A||B|}{1 -18\lambda-o(1)}\right).\]
Also, $|A||B|\ge |F|/n$, so if we take $\lambda=\frac 1{1500}$,

\begin{equation}\label{statement2_bound}
\Pr(X>|A||B|/108)\le \exp\left(-\frac{(K-1)^2}{K+1}\cdot\frac{6 \lambda |F|/n}{1 -18\lambda-o(1)}\right) \le \exp\left( -(1-o(1))\frac{|F|}{500 n} \right).
\end{equation}
Since $|F| \ge M n^2/6$, if we take $M=9000$ we have
\[\Pr(X>|A||B|/108) \le \exp\left( - (1-o(1))\cdot3n\right) = o(8^{-n}).\]
It should be easy to substantially improve the estimate of $M$, but we do not do it here, since we have no specific guess as to the best attainable bound.

\end{proof}

\section{Proof of theorem~\ref{empty:cube}}

Fix $s \geq 100 \cdot \sqrt{n \log n}$. We show that whp the Latin square constructed using Proposition \ref{prop:keevash_empty} has no empty cube of side length $s$.

Indeed, in the spirit of the above proof, we give an upper bound on the probability that a fixed triple $A,B,C \subseteq [n]$ with $|A|=|B|=|C|=s$ is an empty box in a Keevash Steiner triple system. The bound that we get is the probability that statement (\ref{f_safe}) fails plus the probability that $A\times B\times C$ is empty given that statement (\ref{f_safe}) holds, which is at most $2 e^{-(1+o(1))|F|/500 n} \leq 2 e^{-s^3/3000 n}$. Applying the union bound, the probability that there is such a box is at most
\[
\binom{n}{s}^3\cdot 2 e^{-s^3/3000 n} \leq 2 \exp\left(s (3 \log n - s^2/3000 n) \right).
\]
This tends to zero for $s \geq 100 \cdot \sqrt{n \log n}$.

\section{open questions and concluding remarks}\label{section:open}

Let us recall some of the open questions and aims raised above. There are some obvious implications among them, as the reader can easily see.
\begin{enumerate}
\item
Can one find explicit constructions of high-dimensional permutations with good discrepancy properties? Kedlaya's theorem, mentioned above, suggests that substantial new ideas will be required to accomplish this.
\item
Prove that $\varepsilon(L)=O(n^2)$ for {\em almost every} order-$n$ Latin square.
\item
Prove that for every $d\ge 3$ there exist $d$-dimensional permutations $X$ with
$\varepsilon(X)=O(n^2)$.
\item
Prove that $\varepsilon(X)=O(n^2)$ for almost every $d$-dimensional permutation.
\item
Prove the discrepancy conjecture~\ref{conj:high_dim_EML} for Latin squares.
\item
Prove the analogous discrepancy conjecture in all dimensions $d\ge 2$.
\item
Do there exist $d$-dimensional permutations of order-$n$ in which every empty cube has side $\tilde O(n^{1/d})$? Here $\tilde O$ refers to an unspecified polylog term, but perhaps this is even true with $O(n^{1/d})$, which would clearly be tight.
\end{enumerate}

We are presently unable to extend Theorem \ref{thm:typical} to dimensions $d\geq 3$, since the available bounds on the number of $d$-dimensional permutations are not tight enough. Note that this would require very accurate estimates, which seem out of reach with current methods. In this context we recall our conjectured lower bound~\cite{LiLu14} of $\left((1+o(1))\frac{n}{e^d}\right)^{n^d}$. It is conceivable that the machinery in \cite{Ke14} may be useful in this pursuit.

The prospect of extending Theorem \ref{thm:upper_bound} to higher dimensions seems more hopeful. In our proof we show that Keevash's triple systems contain no large empty boxes. This yields this property for the Latin squares representing these Steiner systems. The proof of Proposition \ref{prop:keevash_empty} goes through for $(n,d+1,d)$-Steiner systems in general, and for $d=3$ it is even possible to associate  $3$-dimensional permutations to such Steiner systems. Namely, to an $(n,4,3)$-Steiner system $X$ we associate the $3$-dimensional permutations $A$ given by $A(i,j,k,l)=1$ if $\{i,j,k,l\}\in X$, and $A(i,i,j,j)=1$ for all $i,j \in [n]$. Therefore Theorem \ref{thm:upper_bound} holds in dimension $3$ as well. However, in dimensions $d\geq4$ there seems to be no obvious way of associating Steiner systems with permutations, and so a different approach is needed. It is natural to try and adapt Keevash's method to the construction of high-dimensional permutations, i.e., to analyze the random greedy algorithm in this setting.

We note that the Latin squares constructed in the proof of Theorem \ref{thm:upper_bound} have a large discrepancy, due to \textit{overly dense} boxes that they contain. Keevash's construction associates each vertex $v \in V$ with an element $a_v \in \mathbb{F}_{2^a}$, where $2^{a-2}\leq n\leq 2^{a-1}$. He then considers triples, $x,y,z \in V$ such that $a_x+a_y+a_z=0$ in $\mathbb{F}_{2^a}$. Such a triple that remains legal at the end of the greedy process, gets added to the Steiner triple system. But the additive group of $\mathbb{F}_{2^a}$ has many subgroups. If we take $X=Y=Z$ to be the members of a subgroup, we obtain a collection of vertices with many triples. From the perspective of Latin squares, this is an overly dense box.


It would be interesting to find an explicit construction of Latin squares without large empty boxes. However, most of the known explicit constructions of Latin squares come from groups, but Kedlaya's theorem implies that the multiplication tables of groups always have large empty boxes, which indicates that new ideas are needed here.

Our main conjecture can be viewed as a special case of a much broader problem, that we state in terms of Latin squares, but extensions to permutations of arbitrary dimensions suggest themselves as well. Consider an order-$n$ Latin square $A$, a subset $S\subseteq [n]$ and an index $3\ge t \ge 1$, say $t=2$. The corresponding {\em section} of $A$ is a bipartite graph $\Sigma_{S,t}=(U,V,E)$ on the vertex set $U \cup V$, where $U=V=[n]$ and $ij\in E$ iff there is an $x\in S$ for which $A(i,x,j)=1$. We call this a $k$-section where $|S|=k$. Now pick a parameter of interest $f=f(G)$ that is defined for $k$-regular bipartite graphs $G$ each part of which has $n$ vertices and let $F(k,n)$ be the optimum of $f$ over all such graphs.

{\bf Problem}: Do there exist Latin squares such that $f(\Sigma_{S,t})=(1+o(1))F(k,n)$ for every $k$-section of $A$? For which graph parameters does this hold for almost every Latin square? For the function $f(G) = \max_{A \subset U, B \subset V} |E(A,B)-\frac{k}{n}|A||B||$ we recover our discrepancy conjecture for Latin squares. Many other functions and problems suggest themselves, e.g., minimizing $f(G)$, the largest nontrivial eigenvalue of $G$.

\end{document}